\documentclass[10pt]{amsart}
\usepackage{graphicx,amscd,color,amsmath,amsfonts,amssymb,geometry,soul}

\newtheorem{theorem}{Theorem}
\theoremstyle{plain}

\numberwithin{equation}{section}
\input{tcilatex}

\begin{document}
\title[Summability of multilinear forms on classical sequence spaces]{%
Summability of multilinear forms on classical sequence spaces}
\author[Nogueira]{T. Nogueira}
\address[T. Nogueira]{Departamento de Matem\'{a}tica\\
\indent Universidade Federal da Para\'{\i}ba\\
\indent 58.051-900 - Jo\~{a}o Pessoa, Brazil.}
\email{tonykleverson@gmail.com}
\author[Rueda]{P. Rueda}
\address[P. Rueda]{Departamento de An\'{a}lisis Matem\'{a}tico\\
\indent Universidad de Valencia\\
\indent 46100 Burjassot, Valencia.}
\email{pilar.rueda@uv.es}
\thanks{T. Nogueira is supported by Capes}

\begin{abstract}
We present an extension of the Hardy--Littlewood inequality for multilinear
forms. More precisely, let $\mathbb{K}$ be the real or complex scalar field
and $m,k$ be positive integers with $m\geq k\,$ and $n_{1},\dots ,n_{k}$ be
positive integers such that $n_{1}+\cdots +n_{k}=m$.

($a$) If $(r,p)\in (0,\infty )\times \lbrack 2m,\infty ]$ then there is a
constant $D_{m,r,p,k}^{\mathbb{K}}\geq 1$ (not depending on $n$) such that 
\begin{equation*}
\left( \sum_{i_{1},\dots ,i_{k}=1}^{n}\left\vert T\left(
e_{i_{1}}^{n_{1}},\dots ,e_{i_{k}}^{n_{k}}\right) \right\vert ^{r}\right) ^{%
\frac{1}{r}}\leq D_{m,r,p,k}^{\mathbb{K}} \cdot n^{max\left\{ \frac{%
2kp-kpr-pr+2rm}{2pr},0\right\} }\left\Vert T\right\Vert
\end{equation*}%
for all $m$-linear forms $T:\ell _{p}^{n}\times \cdots \times \ell
_{p}^{n}\rightarrow \mathbb{K}$ and all positive integers $n$. Moreover, the
exponent $max\left\{ \frac{2kp-kpr-pr+2rm}{2pr},0\right\} $ is optimal.

($b$) If $(r, p) \in (0, \infty) \times (m, 2m]$ then there is a constant $%
D_{m,r,p, k}^{\mathbb{K}}\geq 1$ (not depending on $n$) such that 
\begin{equation*}
\left( \sum_{i_{1},\dots ,i_{k}=1}^{n }\left\vert T\left(
e_{i_{1}}^{n_{1}},\dots ,e_{i_{k}}^{n_{k}}\right) \right\vert ^{r }\right) ^{%
\frac{1}{r }}\leq D_{m,r,p, k}^{\mathbb{K}} \cdot n^{ max \left\{\frac{%
p-rp+rm}{pr}, 0\right\}}\left\Vert T\right\Vert
\end{equation*}%
for all $m$-linear forms $T:\ell _{p}^{n}\times \cdots \times \ell
_{p}^{n}\rightarrow \mathbb{K}$ and all positive integers $n$. Moreover, the
exponent $max \left\{\frac{p-rp+rm}{pr}, 0\right\}$ is optimal.

The case $k=m$ recovers a recent result due to G. Araujo and D. Pellegrino.
\end{abstract}

\maketitle

%\subjclass[2010]{Primary 46B25, 47H60}

\section{Introduction}

Let $\mathbb{K}$ be $\mathbb{R}$ or $\mathbb{C},$ $m$ be a positive integer
and $p_{1},...,p_{m}\in \lbrack 1,\infty ]$. For $\mathbf{p}:=(p_{1},\ldots
,p_{m})\in \lbrack 1,\infty ]^{m}$, let 
\begin{equation*}
\left\vert \frac{1}{\mathbf{p}}\right\vert :=\frac{1}{p_{1}}+\cdots +\frac{1%
}{p_{m}},
\end{equation*}%
and let us denote $X_{p}:=\ell _{p}$, when $1\leq p<\infty $, and $X_{\infty
}:=c_{0}$. The following problem has been investigated since the $30^{\prime
}s$ and has important applications:

What is the best value of $\rho $ such that there is a constant $C_{\rho ,%
\mathbf{p}}^{\mathbb{K}}$ such that 
\begin{equation}
\left( \sum\limits_{j_{1},...,j_{m}=1}^{n}\left\vert T\left(
e_{j_{1}},...,e_{j_{m}}\right) \right\vert ^{\rho }\right) ^{\frac{1}{\rho }%
}\leq C_{\rho ,\mathbf{p}}^{\mathbb{K}}\Vert T\Vert  \label{desig}
\end{equation}%
for all continuous $m$-linear forms $T:X_{p_{1}}\times \cdots \times
X_{p_{m}}\rightarrow \mathbb{K}$ and all positive integers $n$?

The answer is divided in some cases; for instance:

\begin{itemize}
\item[(1)] $\rho =\dfrac{2m}{m+1}$, when $\mathbf{p}=(\infty ,...,\infty )$
\quad (Bohnenblust--Hille, \cite{bh});

\item[(2)] $\rho =\dfrac{2m}{m+1-2\left\vert \frac{1}{\mathbf{p}}\right\vert 
}$, when $\left\vert \frac{1}{\mathbf{p}}\right\vert \leq \frac{1}{2}$ \quad
(Hardy--Littlewood \cite{hardy} and Praciano-Pereira \cite{pra});

\item[(3)] $\rho =\dfrac{1}{1-\left\vert \frac{1}{\mathbf{p}}\right\vert }$,
when $\frac{1}{2}\leq \left\vert \frac{1}{\mathbf{p}}\right\vert <1$ \quad
(Hardy--Littlewood \cite{hardy} and Dimant--Sevilla-Peris \cite{dimant});

\item[(4)] $\rho =1$, when $\mathbf{p}=(\infty ,...,\infty )$ and $%
j_{1}=\cdots =j_{m}=j$ \quad (Aron and Globevnik \cite{agg});

\item[(5)] $\rho =\dfrac{1}{1-\left\vert \frac{1}{\mathbf{p}}\right\vert }$,
when $\left\vert \frac{1}{\mathbf{p}}\right\vert <1$ and $j_{1}=\cdots
=j_{m}=j$ \quad (Zalduendo \cite{zal}).
\end{itemize}

These results were successfully unified in a unique inequality in \cite%
{BHHLblocos}, thanks to consider repeated indexes $j_{k}$ in the summands.
Let $n_{1},\ldots ,n_{k}$ be positive integers and $n_{1}+\cdots +n_{k}=m$,
and let us denote by $\left( e_{i_{1}}^{n_{1}},\ldots
,e_{i_{k}}^{n_{k}}\right) $ the $m$-tuple 
\begin{equation*}
(e_{i_{1}},\overset{\text{{\tiny $n_{1}$\thinspace times}}}{\ldots }%
,e_{i_{1}},\ldots ,e_{i_{k}},\overset{\text{{\tiny $n_{k}$\thinspace times}}}%
{\ldots },e_{i_{k}}).
\end{equation*}%
In \cite{BHHLblocos} the following result is proved:

\begin{theorem}[Albuquerque et al. \protect\cite{BHHLblocos}]
\label{thmain2} Let $m\geq k\geq 1,\,m<p\leq \infty $ and let $n_{1},\dots
,n_{k}\geq 1$ be such that $n_{1}+\cdots +n_{k}=m$. Then, for every
continuous $m$--linear form $T:X_{p}\times \cdots \times X_ {p}\rightarrow 
\mathbb{K}$, there is a constant $M(k,m,p,\mathbb{K})\geq 1$ such that 
\begin{equation}  \label{teo1}
\left( \sum_{i_{1},\dots ,i_{k}=1}^{\infty }\left\vert T\left(
e_{i_{1}}^{n_{1}},\dots ,e_{i_{k}}^{n_{k}}\right) \right\vert ^{\rho
}\right) ^{\frac{1}{\rho }}\leq M(k,m,p,\mathbb{K})\left\Vert T\right\Vert ,
\end{equation}
with 
\begin{equation}  \label{expo1}
\rho =\frac{p}{p-m}\text{ for }m<p\leq 2m\text{ and }M(k,m,p,\mathbb{K})\leq
C_{k,p}^{\mathbb{K}}
\end{equation}
and 
\begin{equation}  \label{expo2}
\rho =\frac{2kp}{kp+p-2m}\text{ for }p\geq 2m\text{ and }M(k,m,p,\mathbb{K}
)\leq D_{k,p}^{\mathbb{K}}.
\end{equation}
Moreover, in both cases, the exponent $\rho $ is optimal.
\end{theorem}

The optimality of the exponent in \eqref{teo1}, implies that no constant
independent of $n$ can be found for all $m$-linear forms when a smaller
exponent $r$ is considered. Our objective is to show that, even for smaller
exponents $r$, the value of the left hand sum increases in $n$ under
control, with an explicit dependence on a power factor of $n$. We give
exactly the optimal exponent for $n$. Some previous incursions to this new
approach have been done in \cite[Corollary 5.20]{BoMiPe}. However, it is in 
\cite{archiv} where this subject has been first explored in its own; there,
Hardy-Littlewood type inequalities have been considered and that paper has
been the trigger of our work. More recently, in \cite{GaMaMu} inequalities
involving homogeneous polynomials are studied and the asymptotic behavior of
the constants whenever the number of variables tends to infinity is
established.

This paper is a natural continuation of \cite{BHHLblocos} and, in some sense
it is also related to the notion of index of summability introduced by Maia,
Pellegrino and Santos \cite{maia}, which essentially investigates what
dependence on $n$ emerges when we perturb some well known inequalities (for
instance in (\ref{desig})).

The main result of the present paper is the following:

\begin{theorem}
\label{8j}Let $m,k$ be positive integers, $m\geq k,\,$and let $n_{1},\dots
,n_{k}$ be positive integers such that $n_{1}+\cdots +n_{k}=m$.

(a) If $(r,p)\in (0,\infty )\times \lbrack 2m,\infty ]$ then there is a
constant $D_{m,r,p,k}^{\mathbb{K}}\geq 1$ (not depending on $n$) such that 
\begin{equation*}
\left( \sum_{i_{1},\dots ,i_{k}=1}^{n}\left\vert T\left(
e_{i_{1}}^{n_{1}},\dots ,e_{i_{k}}^{n_{k}}\right) \right\vert ^{r}\right) ^{%
\frac{1}{r}}\leq D_{m,r,p,k}^{\mathbb{K}} \cdot n^{max\left\{ \frac{%
2kp-kpr-pr+2rm}{2pr},0\right\} }\left\Vert T\right\Vert
\end{equation*}%
for all $m$-linear forms $T:X_{p}\times \cdots \times X_{p}\rightarrow 
\mathbb{K}$ and all positive integers $n$. Moreover, the exponent $%
max\left\{ \frac{2kp-kpr-pr+2rm}{2pr},0\right\} $ is optimal.

(b) If $(r, p) \in (0, \infty) \times (m, 2m]$ then there is a constant $%
D_{m,r,p, k}^{\mathbb{K}}\geq 1$ (not depending on $n$) such that 
\begin{equation*}
\left( \sum_{i_{1},\dots ,i_{k}=1}^{n }\left\vert T\left(
e_{i_{1}}^{n_{1}},\dots ,e_{i_{k}}^{n_{k}}\right) \right\vert ^{r }\right) ^{%
\frac{1}{r }}\leq D_{m,r,p, k}^{\mathbb{K}} \cdot n^{ max \left\{\frac{%
p-rp+rm}{pr}, 0\right\}}\left\Vert T\right\Vert
\end{equation*}%
for all $m$-linear forms $T:X_{p}\times \cdots \times X _{p}\rightarrow 
\mathbb{K}$ and all positive integers $n$. Moreover, the exponent $max
\left\{\frac{p-rp+rm}{pr}, 0\right\}$ is optimal.
\end{theorem}

\bigskip

\section{The proof}

Let $E_{1},...,E_{m}$ be Banach spaces. The product $\hat{\otimes}_{j\in
\{1,\ldots ,m\}}^{\pi }E_{j}=E_{1}\hat{\otimes}^{\pi }\cdots \hat{\otimes}%
^{\pi }E_{m}$ denotes the $m$-fold completed projective tensor product of $%
E_{1},\ldots ,E_{m}$. The tensor $x_{1}\otimes \cdots \otimes x_{m}$ will be
denoted by $\otimes _{j\in \{1,\ldots ,m\}}x_{j}$, and $\otimes _{m}x$ shall
denote the tensor $x\otimes \cdots \otimes x$.

%\textcolor{blue}{We can assume $k\geq 2$. The case $k=1$ corresponds with............}
Define $\frac{1}{r_{j}}=\frac{n_{j}}{p}$ and note that $\frac{1}{r_{j}}<1$
for all $j=1,...,k$ (because $p>m$). Let $D_{r_{j}}\subset X_{p}\hat{\otimes}%
^{\pi }\cdots \hat{\otimes}^{\pi }X_{p}$ ($n_{j}$ times) be the vector space
generated by the tensors $\otimes _{n_{j}}e_{i}$ and consider the isometric
isomorphism (see \cite{ArFa} and \cite{BHHLblocos}) $u_{j}:X_{r_{j}}%
\rightarrow \overline{D_{r_{j}}}$ defined by 
\begin{equation*}
u_{j}\left( \sum_{i=1}^{\infty }a_{i}e_{i}\right) =\sum_{i=1}^{\infty
}a_{i}\,\otimes _{n_{j}}\!\! e_{i}.
\end{equation*}

For any continuous $m$-linear form $T:X_{p}\times \cdots \times
X_{p}\rightarrow \mathbb{K}$, consider its $k$-linearization $\widehat{T}:$ $%
\otimes _{n_{1}}^{\pi }X_{p}\times \cdots \times \otimes _{n_{k}}^{\pi
}X_{p} $, that is, $\widehat T$ is the unique $k$-linear form such that $%
\widehat T(x_1^1\otimes \cdots\otimes x_{n_1}^1,\ldots,x_1^k\otimes
\cdots\otimes
x_{n_k}^k)=T(x_1^1,\ldots,x_{n_1}^1,\ldots,x_1^k,\ldots,x_{n_k}^k)$ for all $%
x_j^i\in X_p$, $1\leq j\leq n_i$, $1\leq i\leq k$ (for further details we
refer to \cite{BHHLblocos}), and let $S:X_{r_{1}}\times \cdots \times
X_{r_{k}}\rightarrow \mathbb{K}$ be given by 
\begin{equation*}
S(w_{1},...,w_{k}):=\widehat{T}(u_{1}(w_{1}),...,u_{k}(w_{k})).
\end{equation*}%
Then 
\begin{equation}
\begin{array}{cl}
\left( \displaystyle\sum_{i_{1},\dots ,i_{k}=1}^{n}\left\vert T\left(
e_{i_{1}}^{n_{1}},\dots ,e_{i_{k}}^{n_{k}}\right) \right\vert ^{r}\right) ^{%
\frac{1}{r}} & =\left( \displaystyle\sum_{i_{1},\dots
,i_{k}=1}^{n}\left\vert \widehat{T}\left( u_{1}(e_{i_{1}}),\dots
,u_{k}(e_{i_{k}})\right) \right\vert ^{r}\right) ^{\frac{1}{r}} \\ 
& =\left( \displaystyle\sum_{i_{1},\dots ,i_{k}=1}^{n}\left\vert S\left(
e_{i_{1}},\dots ,e_{i_{k}}\right) \right\vert ^{r}\right) ^{\frac{1}{r}}.%
\end{array}
\label{igual}
\end{equation}

\textbf{Proof of (a).} Let us first suppose that $\left( r,p\right) \in
\left( 0,\frac{2kp}{kp+p-2m}\right] \times \left[ 2m,\infty \right] $. Using
the H\"{o}lder inequality and Theorem \ref{thmain2} we have 
\begin{align*}
& \left( \sum_{i_{1},\dots ,i_{k}=1}^{n}\left\vert T\left(
e_{i_{1}}^{n_{1}},\dots ,e_{i_{k}}^{n_{k}}\right) \right\vert ^{r}\right) ^{%
\frac{1}{r}} \\
& \leq \left( \sum_{i_{1},\dots ,i_{k}=1}^{n}\left\vert T\left(
e_{i_{1}}^{n_{1}},\dots ,e_{i_{k}}^{n_{k}}\right) \right\vert ^{\frac{2kp}{%
kp+p-2m}}\right) ^{\frac{kp+p-2m}{2kp}}\cdot \left( \sum_{i_{1},\ldots
,i_{k}=1}^{n}\left\vert 1\right\vert ^{\frac{2kpr}{2kp-rkp-rp+2mr}}\right) ^{%
\frac{2kp-rkp-rp+2mr}{2kpr}} \\
& =\left( \sum_{i_{1},\dots ,i_{k}=1}^{n}\left\vert T\left(
e_{i_{1}}^{n_{1}},\dots ,e_{i_{k}}^{n_{k}}\right) \right\vert ^{\frac{2kp}{%
kp+p-2m}}\right) ^{\frac{kp+p-2m}{2kp}}\cdot (n^{k})^{\frac{2kp-rkp-rp+2mr}{%
2kpr}} \\
& \leq M(k,m,p,\mathbb{K})\left\Vert T\right\Vert \cdot (n^{k})^{\frac{%
2kp-rkp-rp+2mr}{2kpr}} \\
& =M(k,m,p,\mathbb{K})\Vert T\Vert \cdot n^{\frac{2kp-rkp-rp+2mr}{2pr}}.
\end{align*}

On the other hand, if $\left( r,p\right) \in \left[ \frac{2kp}{kp+p-2m}%
,\infty \right] \times \left[ 2m,\infty \right] $, we have

\begin{align*}
& \left( \sum_{i_{1},\dots ,i_{k}=1}^{n}\left\vert T\left(
e_{i_{1}}^{n_{1}},\dots ,e_{i_{k}}^{n_{k}}\right) \right\vert ^{r}\right) ^{%
\frac{1}{r}} \\
& \leq \left( \sum_{i_{1},\dots ,i_{k}=1}^{\infty }\left\vert T\left(
e_{i_{1}}^{n_{1}},\dots ,e_{i_{k}}^{n_{k}}\right) \right\vert ^{\frac{2kp}{%
kp+p-2m}}\right) ^{\frac{kp+p-2m}{2kp}} \\
& \leq M(k,m,p,\mathbb{K})\left\Vert T\right\Vert \\
& =M(k,m,p,\mathbb{K})\Vert T\Vert \cdot n^{max\left\{ \frac{2kp-rkp-rp+2mr}{%
2pr},0\right\} }
\end{align*}%
and, of course, in this case the exponent $max\left\{ \frac{2kp-rkp-rp+2mr}{%
2pr},0\right\} $ is optimal.

It remains to prove the optimality of the exponent in the case $\left(
r,p\right) \in \left( 0,\frac{2kp}{kp+p-2m}\right] \times \left[ 2m,\infty %
\right] .$ We shall use a technique used in the main result of \cite%
{BHHLblocos}. Suppose that $\lambda \geq 0$ is the smallest exponent
satisfying 
\begin{equation}
\left( \sum_{i_{1},\dots ,i_{k}=1}^{n}\left\vert T\left(
e_{i_{1}}^{n_{1}},\dots ,e_{i_{k}}^{n_{k}}\right) \right\vert ^{r}\right) ^{%
\frac{1}{r}}\leq D_{m,r,p,k}^{\mathbb{K}}\cdot n^{\lambda }\left\Vert
T\right\Vert  \label{otima}
\end{equation}%
for all continuous $m$-linear forms $T:X_{p}\times \cdots \times
X_{p}\rightarrow \mathbb{K}.$ Let us show that $\lambda =max\left\{ \frac{%
2mr+2kp-kpr-pr}{2pr},0\right\} $. Let $A:X_{r_{1}}\times \cdots \times
X_{r_{k}}\rightarrow \mathbb{K}$ be a continuous $k$-linear form. For each $%
i=1,...,k$ we know that $\overline{D}_{r_{i}}$ is complemented into $\hat{%
\otimes}_{j\in \{1,\ldots ,m\}}^{\pi }X_{p}$, and consider the canonical
projection $d_{r_{i}}:$ $\hat{\otimes}_{j\in \{1,\ldots ,m\}}^{\pi
}X_{p}\rightarrow $ $\overline{D}_{r_{i}}$ (see \cite{ArFa} for details).
Defining the $m$-linear form $T_{A}:X_{p}\times \cdots \times
X_{p}\rightarrow \mathbb{K}$ by 
\begin{equation*}
T_{A}(x_{1}^{(1)},\ldots ,x_{n_{1}}^{(1)},\ldots ,x_{1}^{(k)},\ldots
,x_{n_{k}}^{(k)}):=A(u_{r_{1}}^{-1}\circ d_{r_{1}}(x_{1}^{(1)}\otimes \cdots
\otimes x_{n_{1}}^{(1)}),\ldots ,u_{r_{k}}^{-1}\circ
d_{r_{k}}(x_{1}^{(k)}\otimes \cdots \otimes x_{n_{k}}^{(k)})),
\end{equation*}%
we have 
\begin{eqnarray}
T_{A}(e_{i_{1}}^{n_{1}},\ldots ,e_{i_{k}}^{n_{k}}) &=&A(u_{r_{1}}^{-1}\circ
d_{r_{1}}(\otimes _{n_{1}}e_{i_{1}}),\ldots ,u_{r_{k}}^{-1}\circ
d_{r_{k}}(\otimes _{n_{k}}e_{i_{k}}))  \label{qqqq} \\
&=&A(u_{r_{1}}^{-1}(\otimes _{n_{1}}e_{i_{1}}),\ldots
,u_{r_{k}}^{-1}(\otimes _{n_{k}}e_{i_{k}}))=A(e_{i_{1}},\ldots ,e_{i_{k}}). 
\notag
\end{eqnarray}%
By (\ref{qqqq}) and (\ref{otima}) applied to $T_A$, and using that $%
\|T_A\|\leq \|A\|$, we obtain 
\begin{equation*}
\left( \sum_{i_{1},\dots ,i_{k}=1}^{n}\left\vert A(e_{i_{1}},\ldots
,e_{i_{k}})\right\vert ^{r}\right) ^{\frac{1}{r}}\leq D_{m,r,p,k}^{\mathbb{K}%
}\cdot n^{\lambda }\left\Vert A\right\Vert .
\end{equation*}%
Since $A$ is $k$-linear, and%
\begin{equation*}
\frac{1}{r_{j}}=\frac{n_{j}}{p}\leq \frac{m}{2m}=\frac{1}{2},
\end{equation*}%
from the Kahane--Salem--Zygmund inequality (see \cite[Lemma 6.1]{alb} for
details), there is a constant $C_{k}>0$ such that 
\begin{equation*}
n^{\frac{k}{r}}\leq C_{k}M_{k,r,r_{1},...,r_{k}}^{\mathbb{K}}n^{\lambda}n^{%
\frac{k+1}{2}-(\frac{1}{r_{1}}+\cdots +\frac{1}{r_{k}})}.
\end{equation*}%
Making $n\rightarrow \infty $, we have 
\begin{equation*}
\lambda \geq \frac{1}{r_{1}}+\cdots +\frac{1}{r_{k}}+\frac{2k-kr-r}{2r}.
\end{equation*}%
Since $\frac{1}{r_{1}}+\cdots +\frac{1}{r_{k}}=\frac{m}{p}$, we have 
\begin{equation*}
\lambda \geq \frac{m}{p}+\frac{2k-kr-r}{2r}=max\left\{ \frac{2mr+2kp-kpr-pr}{%
2pr},0\right\} .
\end{equation*}

\textbf{Proof of (b).} Since 
\begin{equation*}
\left( \sum_{i_{1},\dots ,i_{k}=1}^{n}\left\vert T\left(
e_{i_{1}}^{n_{1}},\dots ,e_{i_{k}}^{n_{k}}\right) \right\vert ^{r}\right) ^{%
\frac{1}{r}}\leq \left( \sum_{j_{1},\dots ,j_{m}=1}^{n}\left\vert T\left(
e_{j_{1}},\dots ,e_{j_{m}}\right) \right\vert ^{r}\right) ^{\frac{1}{r}},
\end{equation*}%
by \cite[Theorem 1.1(b)]{archiv} we have 
\begin{equation*}
\left( \sum_{i_{1},\dots ,i_{k}=1}^{n}\left\vert T\left(
e_{i_{1}}^{n_{1}},\dots ,e_{i_{k}}^{n_{k}}\right) \right\vert ^{r}\right) ^{%
\frac{1}{r}}\leq D_{m,r,p,k}^{\mathbb{K}}\cdot n^{max\left\{ \frac{p-rp+rm}{%
pr},0\right\} }\left\Vert T\right\Vert .
\end{equation*}

Let us prove the optimality of the exponent. If 
\begin{equation*}
\frac{mr+p-pr}{pr} \leq 0
\end{equation*}%
the optimality of the exponent $max\{(mr+p-pr)/pr,0\}$ is immediate.

Suppose that the inequality holds for a certain exponent $s\geq 0;$ thus 
\begin{equation}
\left( \sum_{i_{1},\dots ,i_{k}=1}^{n}\left\vert T\left(
e_{i_{1}}^{n_{1}},\dots ,e_{i_{k}}^{n_{k}}\right) \right\vert ^{r}\right) ^{%
\frac{1}{r}}\leq D_{m,r,p,k}^{\mathbb{K}}\cdot n^{s}\left\Vert T\right\Vert .
\label{coms1}
\end{equation}

As in the previous case, for each continuous $k$-linear form $%
A:X_{r_{1}}\times \cdots \times X_{r_{k}}\rightarrow \mathbb{K}$, with $%
r_{j}=\frac{p}{n_{j}},$ and for all $j=1,...,k$, there is a continuous $m$%
-linear form $T_{A}:X_{p}\times \cdots \times X_{p}\rightarrow \mathbb{K}$
such that 
\begin{equation}
T_{A}(e_{i_{1}}^{n_{1}},\ldots ,e_{i_{k}}^{n_{k}})=A(e_{i_{1}},\ldots
,e_{i_{k}})  \label{rsrs}
\end{equation}%
and $\left\Vert T_{A}\right\Vert \leq \left\Vert A\right\Vert $. By (\ref%
{rsrs}) and (\ref{coms1}) applied to $T_A$ we obtain 
\begin{equation}
\left( \sum_{i_{1},\dots ,i_{k}=1}^{n}\left\vert A(e_{i_{1}},\ldots
,e_{i_{k}})\right\vert ^{r}\right) ^{\frac{1}{r}}\leq D_{m,r,p,k}^{\mathbb{K}%
}\cdot n^{s}\left\Vert T_{A}\right\Vert \leq D_{m,r,p,k}^{\mathbb{K}}\cdot
n^{s}\left\Vert A\right\Vert .  \label{y55}
\end{equation}%
Define the $k$-linear form $S:X_{r_{1}}\times \cdots \times
X_{r_{k}}\rightarrow \mathbb{K}$ by 
\begin{equation*}
S(x^{(1)},...,x^{(k)})=\sum_{j=1}^{n}x_{j}^{(1)}\cdots x_{j}^{(k)},
\end{equation*}%
and notice that by the H\"{o}lder inequality we have 
\begin{equation*}
\Vert S\Vert \leq n^{1-\left( \frac{1}{r_{1}}+\cdots +\frac{1}{r_{k}}\right)
}=n^{1-\frac{m}{p}}.
\end{equation*}%
Therefore, plugging $S$ into (\ref{y55}) we get 
\begin{equation*}
n^{\frac{1}{r}}\leq D_{m,r,p,k}^{\mathbb{K}}n^{s}n^{1-\frac{m}{p}}
\end{equation*}%
and we easily conclude that 
\begin{equation*}
s\geq \frac{p-rp+rm}{pr}.
\end{equation*}

\end{document}